\documentclass{ifacconf}

\usepackage{graphicx,picture,calc}
\usepackage{natbib}
\usepackage{amsmath}
\usepackage{amssymb}
\usepackage{mathtools}
\usepackage{enumitem}
\usepackage{tikz}
\usepackage{booktabs}
\usepackage{array}
\usepackage{epstopdf}

\setlist[enumerate]{label={\textnormal{(\alph*)}}, ref={(\alph*)}, leftmargin=*, nosep}

\DeclareMathOperator{\Real}{Re} 
\DeclareMathOperator{\Imag}{Im} 

\newcommand{\suchthat}{\ifnum\currentgrouptype=16 \mathrel{}\middle|\mathrel{}\else\mid\fi}

\graphicspath{{./Figures/}}

\begin{document}
\usetikzlibrary{arrows}
\setlength{\parskip}{6pt plus 2pt minus 1pt} 

\begin{frontmatter}

\title{Spectral dominance of complex roots for single-delay linear equations\thanksref{footnoteinfo}} 

\thanks[footnoteinfo]{Corresponding author: Guilherme Mazanti (gui\-lher\-me.ma\-zan\-ti@l2s.cen\-trale\-sup\-elec.fr)}

\author[LSS,IPSA]{Guilherme Mazanti} 
\author[LSS,IPSA]{Islam Boussaada} 
\author[LSS]{Silviu-Iulian Niculescu}
\author[CVUT]{Tom\'{a}\v{s} Vyhl\'{\i}dal}

\address[LSS]{University Paris-Saclay, CNRS, CentraleSup\'elec, Laboratoire des Signaux et Syst\`emes (L2S), Inria Saclay, DISCO Team, 3, rue Joliot Curie, 91192, Gif-sur-Yvette, France.\\\{first name.last name\}@l2s.centralesupelec.fr}
\address[IPSA]{Institut Polytechnique des Sciences Avanc\'ees (IPSA)\\63 boulevard de Brandebourg, 94200 Ivry-sur-Seine, France.}
\address[CVUT]{Dept.\ of Instrumentation and Control Eng., and Center of Advanced Aerospace Technology, Faculty of Mechanical Engineering, Czech Technical University in Prague.\\Technick\'a 4, 166 07, Praha 6, Czechia. \\tomas.vyhlidal@fs.cvut.cz}

\begin{abstract}
This paper provides necessary and sufficient conditions for the existence of a pair of complex conjugate roots, each of multiplicity two, in the spectrum of a linear time-invariant single-delay equation of retarded type. This pair of roots is also shown to be always strictly dominant, determining thus the asymptotic behavior of the system. The proof of this result is based on the corresponding result for real roots of multiplicity four, continuous dependence of roots with respect to parameters, and the study of crossing imaginary roots. We also present how this design can be applied to vibration suppression and flexible mode compensation.
\end{abstract}


\end{frontmatter}

\section{Introduction}

In this paper, we consider a linear time-invariant equation with a single delay of the form
\begin{equation}
\label{DiffEqn}
y^{\prime\prime}(t) + a_1 y^\prime(t) + a_0 y(t) + \alpha_1 y^\prime(t - \tau) + \alpha_0 y(t - \tau) = 0,
\end{equation}
where the coefficients $a_1, a_0, \alpha_1, \alpha_0$ are real numbers and the delay $\tau$ is a positive real number. Equations of the form \eqref{DiffEqn} are said to be delayed equations of \emph{retarded} type since the derivative of highest order only appears in the non-delayed term $y^{\prime\prime}(t)$.

Time delays are useful for modeling propagation phenomena, such as of material, energy, or information, with a finite propagation speed, this propagation taking place typically between parts of a complex system. For this reason, equations and systems with time delays have been widely used in several scientific and technological domains in which modeling such propagation phenomena is important, such as in biology, chemistry, economics, physics, or engineering. Due to these applications and the challenging mathematical problems arising in their analysis, time-delay systems have been the subject of much attention by researchers in several fields, in particular since the 1950s and 1960s, such as, for instance, in \cite{Bellman1963Differential, Halanay1966Differential, Pinney1958Ordinary}. We refer to \cite{Diekmann1995Delay, Gopalsamy1992Stability, Gu2003Stability, Hale1993Introduction, Insperger2011Semi, Li2017Frequency, Michiels2014Stability, Stepan1989Retarded} for details on time-delay systems and their applications.

The stability analysis of time-delay systems has attracted much research effort and is an active field (see, e.g., \cite{Abdallah1993Delayed, Chen1995New, Cooke1986Zeroes, Gu2003Stability, Michiels2014Stability, Olgac2002Exact, Sipahi2011Stability}). A usual technique for addressing stability of linear time-invariant systems in the delay-free situation is based on spectral methods and consists in considering the corresponding characteristic polynomial, whose complex roots determine the asymptotic behavior of solutions of the system. This technique also carries over for linear time-invariant systems with delays, whose asymptotic behavior can also be characterized in terms of complex roots of a certain characteristic function (see, e.g., \cite[Proposition~1.13]{Michiels2014Stability}). For \eqref{DiffEqn}, this characteristic function is
\begin{equation}
\label{Delta}
\Delta(s) = s^2 + a_1 s + a_0 + e^{- s \tau}(\alpha_1 s + \alpha_0).
\end{equation}
Similarly to the delay-free case, all solutions of \eqref{DiffEqn} converge exponentially fast to $0$ if and only if $\Real s < 0$ for every $s \in \mathbb C$ such that $\Delta(s) = 0$, and the asymptotic behavior of solutions of \eqref{DiffEqn} is determined by the real number $\gamma_0 = \sup\{\Real s \suchthat s \in \mathbb C,\; \Delta(s) = 0\}$, called the \emph{spectral abscissa} of $\Delta$.

Entire functions such as $\Delta$ that can be written under the form $Q(s) = \sum_{k=1}^\ell p_k(s) e^{\lambda_k s}$ for some polynomials with real coefficients $p_1, \dotsc, p_\ell$ and pairwise distinct real numbers $\lambda_1, \dotsc, \lambda_\ell$ are called \emph{quasipolynomials}. The interest in studying quasipolynomials come from the fact that, when $\lambda_k \leq 0$ for every $k$, they are characteristic equations of linear time-invariant delayed equations.

One usually defines the \emph{degree} of a quasipolynomial $Q$ as above to be $D = \ell + \delta - 1$, where $\delta$ is the sum of the degrees of $p_1, \dotsc, p_\ell$ (see, e.g., \cite{Wielonsky2001Rolle, Berenstein1995Complex}). In particular, the degree of $\Delta$ in \eqref{Delta} is $D = 2 + 3 - 1 = 4$. Contrarily to the case of polynomials, the degree of a quasipolynomial does not determine the number of roots of the quasipolynomial, which is infinite except in trivial cases. However, similarly to polynomials, the degree does have a link with multiplicities of roots. More precisely, a classical result on quasipolynomials provided in \cite[Problem 206.2]{Polya1998Problems}, known as the \emph{P\'{o}lya--Szeg\H{o} bound}, implies that, given a quasipolynomial $Q$ of degree $D \geq 0$, the multiplicity of any root of $Q$ does not exceed $D$. For the quasipolynomial $\Delta$ from \eqref{Delta}, this means that any of its roots has multiplicity at most $4$. Recent works such as \cite{Boussaada2016Characterizing, Boussaada2016Tracking} have provided characterizations of multiple roots of quasipolynomials using approaches based on Birkhoff and Vandermonde matrices.

When studying the roots of a quasipolynomial in order to analyze the stability of a time-delay system, only the rightmost roots on the complex plane are important for determining the system's asymptotic behavior. These roots are usually called \emph{dominant roots} and can be defined as follows.
\begin{defn}
Let $Q: \mathbb C \to \mathbb C$ and $s_0 \in \mathbb C$.
\begin{enumerate}
\item We say that $s_0$ is a \emph{dominant} (respectively, \emph{strictly dominant}) root of $Q$ if $Q(s_0) = 0$ and, for every $s \in \mathbb C \setminus \{s_0\}$ such that $Q(s) = 0$, one has $\Real s \leq \Real s_0$ (respectively, $\Real s < \Real s_0$).
\item We say that $s_0$ and its complex conjugate $\overline s_0$ are a pair of \emph{dominant} (respectively, \emph{strictly dominant}) roots of $Q$ if $Q(s_0) = Q(\overline s_0) = 0$ and, for every $s \in \mathbb C \setminus \{s_0, \overline s_0\}$ such that $Q(s) = 0$, one has $\Real s \leq \Real s_0$ (respectively, $\Real s < \Real s_0$).
\end{enumerate}
\end{defn}
Dominant roots may not exist in general, but they always exist for functions of the form \eqref{Delta} (see, e.g., \cite[Chapter 1, Lemma 4.1]{Hale1993Introduction}). Exponential stability of \eqref{DiffEqn} is equivalent to the dominant roots of $\Delta$ having negative real part.

It has been observed in several works that real roots of high multiplicity tend to be dominant, a property known as \emph{multiplicity-induced dominance} (MID for short). We refer the reader, for instance, to \cite{Boussaada2018Further}, in which MID was proved for \eqref{Delta} in the case $\alpha_1 = 0$ for a real root of multiplicity $3$ thanks to a suitable factorization of $\Delta$, and to \cite{Boussaada2020Multiplicity}, which considers the case $\alpha_1 \neq 0$ and proves dominance of a real root of multiplicity $4$ using Cauchy's argument principle. MID is also reminiscent of the fact that, for delay-free systems with an affine constraint on their coefficients, the spectral abscissa is minimized on a polynomial with a single root of maximal multiplicity (see \cite{Blondel2012Explicit, Chen1979Output}), with similar properties for some time-delay systems obtained in \cite{Michiels2002Continuous, Ramirez2016Design, Vanbiervliet2008Nonsmooth}. The interest in considering multiple roots does not rely on the multiplicity itself, but rather on its connection with dominance and the corresponding implications for stability analysis and control design.

The main goal of this paper is to investigate whether MID holds for $\Delta$ when assigning a pair of complex conjugate roots instead of a real root. Designing a system to have a pair of dominant complex conjugate roots may have several practical interests, as highlighted in \cite{Kure2018Spectral}, in which a robust delayed resonator is designed by assigning double imaginary roots, and as we also illustrate in Section~\ref{SecApplications}. The questions we address in this paper are the following.
\begin{enumerate}[label={\textbf{(Q\arabic*)}}, ref={(Q\arabic*)}, nosep, leftmargin=*]
\item\label{Ques1} Is it possible to choose $a_1, a_0, \alpha_1, \alpha_0 \in \mathbb R$ in such a way that a given complex number $s_0$ and its complex conjugate $\overline s_0$ are roots of multiplicity $2$ of $\Delta$?
\item\label{Ques2} Under the above choice, do $s_0$ and $\overline s_0$ form a pair of (strictly) dominant roots?
\end{enumerate}
Our main result, Theorem~\ref{MainTheo}, in addition to recalling the situation for real root assignment, also provides affirmative answers to both questions. Question~\ref{Ques1} can be addressed in a straightforward manner, whereas the answer to \ref{Ques2} relies on the continuity of the other roots of $\Delta$ with respect to the assigned root and a study of crossing imaginary roots, using techniques similar in spirit to those of \cite{Boussaada2016Tracking}.

The paper is organized as follows: Notations used in the paper are standard. Section~\ref{SecMainResult} provides the statement of our main result, Theorem~\ref{MainTheo}, as well as a sketch of its proof, while Section~\ref{SecExpl} contains illustrative examples.

\section{Main result}
\label{SecMainResult}

The main result of this paper is the following.

\begin{thm}
\label{MainTheo}
Consider the quasipolynomial $\Delta$ given by \eqref{Delta} and let $s_0 \in \mathbb C$, $\sigma_0 = \Real s_0$, and $\theta_0 = \Imag s_0$.
\begin{enumerate}
\item\label{ConditionsReal} Assume that $\theta_0 = 0$. Then $s_0$ is a root of multiplicity $4$ of $\Delta$ if and only if the coefficients $a_0, a_1, \alpha_0, \alpha_1$, the value $\sigma_0$, and the delay $\tau$ satisfy the relations
\begin{subequations}
\label{CoeffsReal}
\begin{align}
a_1 & = -\frac{4}{\tau} - 2 \sigma_0, & \qquad a_0 & = \frac{6}{\tau^2} + \frac{4}{\tau} \sigma_0 + \sigma_0^2, \displaybreak[0] \label{CoeffsRealA} \\
\alpha_1 & = -\frac{2}{\tau} e^{\sigma_0 \tau}, & \alpha_0 & = \frac{2}{\tau} e^{\sigma_0 \tau} \left(\sigma_0 - \frac{3}{\tau}\right). \label{CoeffsRealAlpha}
\end{align}
\end{subequations}
\item\label{ConditionsComplex} Assume that $\theta_0 \neq 0$. Then $s_0$ and $\overline s_0$ are roots of multiplicity $2$ of $\Delta$ if and only if the coefficients $a_0, a_1, \alpha_0, \alpha_1$, the values $\sigma_0$ and $\theta_0$, and the delay $\tau$ satisfy the relations
\begin{subequations}
\label{CoeffsComplex}
\begin{align}
a_1 & = - 2 \sigma_{0} - 2 \theta_{0} \frac{\tau \theta_{0} - \sin{\left (\tau \theta_{0} \right )}\cos{\left (\tau \theta_{0} \right )}}{\tau^{2} \theta_{0}^{2} - \sin^{2}{\left (\tau \theta_{0} \right )}}, \displaybreak[0] \label{CoeffsComplexA1}\\
a_0 & = \sigma_{0}^{2} + 2 \sigma_{0} \theta_{0}\frac{\tau \theta_{0} - \sin{\left (\tau \theta_{0} \right )}\cos{\left (\tau \theta_{0} \right )}}{\tau^{2} \theta_{0}^{2} - \sin^{2}{\left (\tau \theta_{0} \right )}} \notag \\
 & \hphantom{{} = {}} + \theta_{0}^{2} \frac{\tau^{2} \theta_{0}^{2} + \sin^{2}{\left (\tau \theta_{0} \right )}}{\tau^{2} \theta_{0}^{2} - \sin^{2}{\left (\tau \theta_{0} \right )}}, \displaybreak[0] \label{CoeffsComplexA0} \\
\alpha_1 & = 2 \theta_{0} e^{\sigma_{0} \tau} \frac{\tau \theta_{0} \cos{\left (\tau \theta_{0} \right )} - \sin{\left (\tau \theta_{0} \right )}}{\tau^{2} \theta_{0}^{2} - \sin^{2}{\left (\tau \theta_{0} \right )}}, \displaybreak[0] \label{CoeffsComplexAlpha1} \\
\alpha_0 & = 2 \theta_{0} e^{\sigma_{0} \tau} \Biggl(\sigma_{0} \frac{\sin{\left (\tau \theta_{0} \right )} - \tau \theta_{0} \cos{\left (\tau \theta_{0} \right )}}{\tau^{2} \theta_{0}^{2} - \sin^{2}{\left (\tau \theta_{0} \right )}} \notag \\
& \hphantom{ = 2 \theta_{0} e^{\sigma_{0} \tau} \Biggl(} - \frac{\tau \theta_{0}^{2} \sin{\left (\tau \theta_{0} \right )}}{\tau^{2} \theta_{0}^{2} - \sin^{2}{\left (\tau \theta_{0} \right )}}\Biggr). \label{CoeffsComplexAlpha0}
\end{align}
\end{subequations}
\item\label{DominanceReal} If \eqref{CoeffsReal} is satisfied, then $s_0$ is a strictly dominant root of $\Delta$.
\item\label{DominanceComplex} If \eqref{CoeffsComplex} is satisfied, then $s_0$ and $\overline s_0$ are a pair of strictly dominant roots of $\Delta$.
\end{enumerate}
\end{thm}

\begin{rem}
The expressions of $a_1, a_0, \alpha_1, \alpha_0$ in \eqref{CoeffsReal} and \eqref{CoeffsComplex} are singular with respect to $\tau$ as $\tau \to 0$. If one is interested in studying the behavior of the roots of $\Delta$ as $\tau \to 0$ when \eqref{CoeffsReal} or \eqref{CoeffsComplex} is satisfied, one may consider instead the quasipolynomial $s \mapsto \tau^2 \Delta(s)$, which has the same roots as $\Delta$ but whose coefficients are regular with respect to $\tau$.
\end{rem}

\begin{rem}
\label{RemkLimit}
The expressions of $a_1, a_0, \alpha_1, \alpha_0$ in \eqref{CoeffsComplex} are well-defined for every $\theta_0 \in \mathbb R \setminus \{0\}$ and $\tau > 0$, since $\sin^2(\tau\theta_0) = \tau^2\theta_0^2$ if and only if $\tau\theta_0 = 0$. Moreover, these expressions are even functions of $\theta_0$ --- as one might expect by symmetry since one is placing both roots $s_0$ and $\overline s_0$ --- and they converge to the corresponding expressions in \eqref{CoeffsReal} as $\theta_0 \to 0$.
\end{rem}

Up to a translation and a scaling of the spectrum represented by the change of variables $z = \tau(s - \sigma_0)$, one may reduce to the case $\sigma_0 = 0$ and $\tau = 1$, in which \eqref{CoeffsReal} reduces to $a_1 = -4$, $a_0 = 6$, $\alpha_1 = -2$, $\alpha_0 = -6$, yielding the quasipolynomial
\begin{equation}
\label{DeltaHatReal}
\widehat\Delta_{\mathrm{R}}(z) = z^2 - 4z + 6 - e^{-z} (2 z + 6),
\end{equation}
and \eqref{CoeffsComplex} reduces to
\begin{subequations}
\label{AAndAlphaNormalized}
\begin{align}
a_1 & = -2 \theta_0 \frac{\theta_0 - \sin\theta_0 \cos\theta_0}{\theta_0^2 - \sin^2 \theta_0}, & a_0 & = \theta_0^2 \frac{\theta_0^2 + \sin^2 \theta_0}{\theta_0^2 - \sin^2 \theta_0}, \displaybreak[0] \\
\alpha_1 & = 2 \theta_0 \frac{\theta_0 \cos \theta_0 - \sin\theta_0}{\theta_0^2 - \sin^2 \theta_0}, & \alpha_0 & = - \frac{2 \theta_0^3 \sin \theta_0}{\theta_0^2 - \sin^2 \theta_0},
\end{align}
\end{subequations}
yielding the quasipolynomial
\begin{multline}
\label{DeltaHatComplex}
\widehat\Delta_{\mathrm C}(z; \theta_0) = z^2 -2 \theta_0 \frac{\theta_0 - \sin\theta_0 \cos\theta_0}{\theta_0^2 - \sin^2 \theta_0} z + \theta_0^2 \frac{\theta_0^2 + \sin^2 \theta_0}{\theta_0^2 - \sin^2 \theta_0} \\ + e^{-z}\Biggl(2 \theta_0 \frac{\theta_0 \cos \theta_0 - \sin\theta_0}{\theta_0^2 - \sin^2 \theta_0} z - \frac{2 \theta_0^3 \sin \theta_0}{\theta_0^2 - \sin^2 \theta_0}\Biggr).
\end{multline}
In the sequel of the paper, we use the convention, in accordance with Remark~\ref{RemkLimit}, that $\widehat\Delta_{\mathrm C}(\cdot; 0) = \widehat\Delta_{\mathrm R}(\cdot)$.

We now provide the main ideas for the proof of Theorem~\ref{MainTheo}. The complete proof can be found in an upcoming extended version of this paper.

\begin{pf*}{Sketch of the proof.} 
We consider only the case $s_0 = i \theta_0$ for $\theta_0 \geq 0$, since the general case can be reduced to it by the above change of variables. Assertions \ref{ConditionsReal} and \ref{DominanceReal} have already been proved in \cite{Boussaada2020Multiplicity, MazantiQualitative}. To prove assertion \ref{ConditionsComplex}, we notice that, for real coefficients $a_1, a_0, \alpha_1, \alpha_0$, $s_0 = i \theta_0$ with $\theta_0 > 0$ satisfies $\Delta(s_0) = \Delta^\prime(s_0) = 0$ if and only if \eqref{AAndAlphaNormalized} holds. In this case, one verifies that $\Delta^{\prime\prime}(s_0) \neq 0$, showing that the multiplicity of $s_0$ is indeed $2$.

To prove \ref{DominanceComplex}, note that, by \ref{DominanceReal}, $0$ is a strictly dominant root of $\widehat\Delta_{\mathrm C}(\cdot; 0)$ of multiplicity $4$ and, as $\theta_0$ increases, this root splits into two roots $\pm i \theta_0$ of $\widehat\Delta_{\mathrm C}(\cdot; \theta_0)$ of multiplicity $2$ each. By continuity of the roots of $\widehat\Delta_{\mathrm C}(\cdot; \theta_0)$ with respect to $\theta_0$, the roots $\pm i \theta_0$ will cease to be strictly dominant as $\theta_0$ increases if and only if one root coming from $\infty$ appears at the right half-plane or one root with negative real part crosses the imaginary axis. The proof is completed by arguing by contradiction to show that none of these two cases may occur.
\end{pf*}

\section{Illustrative examples}
\label{SecExpl}

\subsection{Roots of $\widehat\Delta_{\mathrm C}(\cdot; \theta_0)$ as a function of $\theta_0$}
\label{SecRootsFuncTheta}

The quasipolynomial $\widehat\Delta_{\mathrm C}(\cdot; \theta_0)$ from \eqref{DeltaHatComplex} is obtained by applying Theorem~\ref{MainTheo} to $s_0 = i\theta_0$ for some $\theta_0 \in \mathbb R$. Theorem~\ref{MainTheo} guarantees that the multiple roots $\pm i \theta_0$ are strictly dominant, but says nothing about how the roots on the open left half-plane behave. In order to get a grasp on their behavior, we have performed numerical computations of all roots of $\widehat\Delta_{\mathrm C}(\cdot; \theta_0)$ on the region $\{s \in \mathbb C \suchthat -4.75 \leq \Real s \leq 0.25 \text{ and } -25 \leq \Imag s \leq 25\}$ for several values of $\theta_0 \in [0, 8]$. The results are provided in Fig.~\ref{FigRoots}, with different values of $\theta_0$ being represented with different colors. All numerical computations have been performed using Python \texttt{cxroots} package, which implements numerical methods described in \cite{Kravanja2000Computing}.

\begin{figure}[ht]
\centering
\resizebox{\columnwidth}{!}{\input{Figures/roots.pgf}}
\caption{Roots of $\widehat\Delta_{\mathrm C}(\cdot; \theta_0)$ for $\theta_0 \in [0, 8]$, with a detailed view of the region $\{s \in \mathbb C \suchthat -1.80 \leq \Real s \leq -1.55 \text{ and } 10 \leq \Imag s \leq 11\}$.}
\label{FigRoots}
\end{figure}

One observes in Fig.~\ref{FigRoots} the movement of the dominant roots $\pm i \theta_0$ along the imaginary axis. Concerning the other roots, as $\theta_0$ increases, the imaginary parts of the non-real roots represented in the figure increase in absolute value, while the real parts oscillate. Table~\ref{TabFirstRoot} synthesizes this oscillatory behavior for the first pair of non-dominant complex conjugate roots in Fig.~\ref{FigRoots} (in order of increasing absolute value of imaginary part) by presenting the values of $\theta_0$ and the corresponding roots at local extrema of their real part, as well as for the initial and final values $\theta_0 = 0$ and $\theta_0 = 8$ used in the numerical computations. 

\begin{table}[ht]
\centering
\caption{First pair of non-dominant complex conjugate roots in Fig.~\ref{FigRoots} at their initial and final positions and on local extrema of their real parts.}
\label{TabFirstRoot}
\begin{tabular}{ccc}
\toprule
 & $\theta_0$ & Roots \\
\midrule
Initial & $0$ & $-1.731 \pm 10.16 i$ \\
First local maximum & $2.51$ & $-1.586 \pm 10.46i$ \\
First local minimum & $4.59$ & $-2.735 \pm 12.14i$ \\
Second local maximum & $6.19$ & $-1.764 \pm 13.74i$ \\
Second local minimum & $7.83$ & $-2.508 \pm 15.32i$ \\
Final & $8$ & $-2.466 \pm 15.66i$ \\
\bottomrule
\end{tabular}
\end{table}

\begin{figure*}[ht]
\centering
\resizebox{\textwidth}{!}{\input{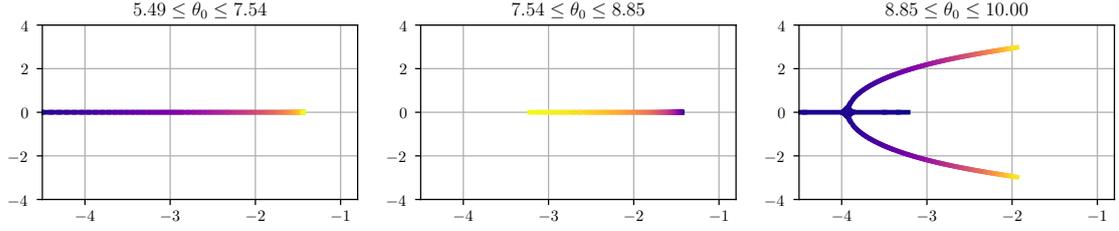}}
\caption{Detailed view of the roots of $\widehat\Delta_{\mathrm C}(\cdot; \theta_0)$ around the real axis for different ranges of $\theta_0$. In each figure, darker colors correspond to smaller values of $\theta_0$.}
\label{FigRealRoots}
\end{figure*}

We also notice in Fig.~\ref{FigRoots} the presence of a real-valued root. Its detailed behavior obtained from numerical computations for $\theta_0 \in [5.49, 10.00]$ is provided in Fig.~\ref{FigRealRoots}, which is split in three different ranges for $\theta_0$ corresponding to different observed behaviors of the root. This root first appears in the domain under consideration for $\theta_0 \approx 5.49$ and moves to the right, reaching a local maximum at $\theta_0 \approx 7.54$, at which point its value is approximately $-1.437$. It then starts moving to the left for $\theta_0 \in [7.54, 8.85]$. At $\theta_0 \approx 8.85$, a second real root appears in the domain under consideration, coming from $-\infty$ and moving to the right, and both roots meet, giving rise, when $\theta_0 \approx 8.88$, to a real root of multiplicity $2$ whose value is approximately $-3.927$. For $\theta_0 \in [8.88, 10]$, these roots become a pair of complex conjugate roots which start moving to the right. As $\theta_0$ increases beyond $10$ (not represented in Fig.~\ref{FigRealRoots}), one observes that this pair of roots oscillates like the other pairs of complex conjugate roots from Fig.~\ref{FigRoots}.

\subsection{Applications: vibration suppression and fle\-xi\-ble mode compensation}
\label{SecApplications}

We provide two engineering applications, with a common requirement for having a double root on the imaginary axis. The first application is {\em active vibration suppression} (AVS) and the second application is {\em flexible mode compensation} (FMC). The common feature of these two methods is that the purely imaginary roots $\pm i \omega$ of \eqref{Delta} are turned to imaginary zeros of the overall system. In AVS, $\omega$ is the frequency of an excitation force, while, for FMC, $\omega$ is the natural frequency of the flexible mode to be compensated. In both cases, the overall system magnitude at frequency $\omega$ is zero. The multiplicity two of the zero then increases the robustness in the vibration suppression or mode compensation. Before explaining these two applications in more detail, let us propose delay values $\tau$. From the practical point of view, an intuitive choice for the delay is given by
\begin{equation}
    \tau_k=\frac{k\pi}{\omega}, \qquad k = 1,2,3,\dotsc,
    \label{eq:tau_k}
\end{equation}
for which \eqref{CoeffsComplex} gives
$a_1 = - \frac{2 \omega}{k\pi}$, 
$a_0 = \omega^2$,
$\alpha_1 = (-1)^k \frac{2 \omega}{k\pi}$, and
$\alpha_0 = 0$. 
Thus, the characteristic function \eqref{Delta} turns to
\begin{equation}
\label{Delta_w}
\Delta_\omega(s) = s^2 + \frac{2}{\tau_k}\left((-1)^k e^{- s \tau_k}-1\right)s + \omega^2.
\end{equation}

\begin{figure}[htb]
\centering
\includegraphics[width=0.45\columnwidth]{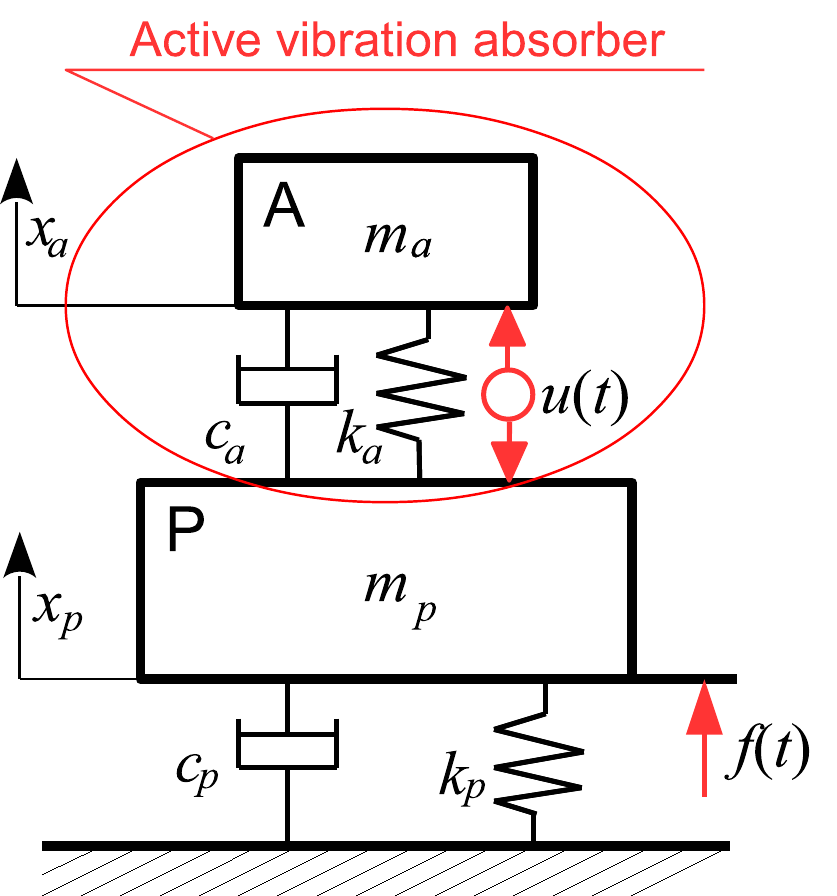}
\caption{Primary Structure (P), with an active vibration absorber (A) to suppress displacement $x_p$ induced by harmonic disturbance force $f(t)$}
\label{fig:intro_1}
\end{figure}
In the AVS application, we adapt the {\em delayed resonator scheme} proposed by \cite{olgac1994novel} with a single root at $\pm i\omega$. Recently the concept has been adjusted by \cite{Kure2018Spectral} with double roots at $\pm i\omega$ in order to enhance the robustness. Let us note that the solution in \cite{Kure2018Spectral} required two time delays. Here, we provide a solution with a single delay. The scheme of the set-up is shown in Fig.~\ref{fig:intro_1}. The system main body is a vibrating platform $P$ excited by a periodic external force $f(t)=F \cos(\omega t)$, $F$ denoting the force amplitude. In order to compensate fully the vibrations, the absorber $A$ is actuated with the active feedback $u(t)$. The absorber dynamics is then given by 
\begin{equation}
x_a^{\prime\prime}\left(t\right)+2\zeta\Omega x_a^\prime \left(t\right)+\Omega^2 x_a \left(t\right)=\frac{1}{m_a} u(t).
\label{eq:Abs}
\end{equation}
where $\zeta,\Omega, m_a$ are the damping, natural frequency and mass of the physical absorber. Introducing the active feedback in the form
\begin{multline}
u(t)=m_a(\Omega^2-\omega^2)x_a(t)+2m_a\left(\zeta\Omega+\frac{1}{\tau_k}\right)x_a^\prime(t)\\
-2m_a\frac{(-1)^k}{\tau_k}x_a^\prime(t-\tau_k),
\label{fdb_DR}
\end{multline}
the characteristic function of the active absorber \eqref{eq:Abs}--\eqref{fdb_DR} is given by \eqref{Delta_w} with a double root at $\pm i\omega$. As demonstrated e.g.\ in \cite{Kure2018Spectral}, the transfer function $f\rightarrow x_p$ is in the form 
\begin{equation}
    G_{x_af}(s)=\frac{\Delta_\omega(s)}{M(s)}
    \label{tf_P_DR}
\end{equation}
where $M(s)$ is a characteristic function of the closed loop system. Therefore, as required, the double roots at $\pm i\omega$ become double zeros of \eqref{tf_P_DR}. This implies that no vibrations are transferred from $f$ to $x_p$ and the platform is fully silenced. 

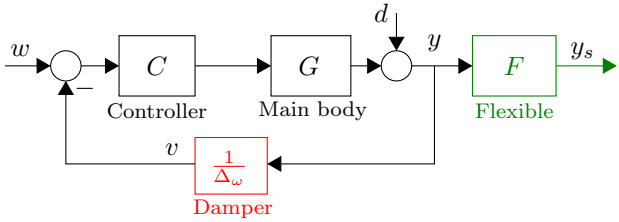
\begin{figure}[ht]
\begin{tikzpicture}
\tikzstyle{myedgestyle} = [-triangle 60]
\node[green!50!black] at (-2.1,1.7) {$F$};
\draw [-triangle 60,green!50!black](-1.55,1.7) -- (-0.75,1.7);
\draw [green!50!black] (-2.65,1.3) rectangle (-1.55,2.1);
\node[red] at (-5.8,-0.2) {\small Damper};
\node[green!50!black] at (-2.1,1.1) {\small Flexible};
\node at (-6.8,1.1) {\small Controller};
\node at (-4.75,1.1) {\small Main body};
\node at (-1.2,1.9) {$y_s$};
\draw [-triangle 60] (-8,1.7) ellipse (0.2 and 0.2);
\draw[red] [-triangle 60] (-6.3,0.75) rectangle (-5.35,0);
\draw [-triangle 60] (-7.3,2.1) rectangle (-6.3,1.3);
\draw [-triangle 60](-3.15,1.7) -- (-3.15,0.4) -- (-5.35,0.4);
\draw [-triangle 60](-6.3,0.4) -- (-8,0.4) -- (-8,1.5);
\draw [-triangle 60](-8.8,1.7) -- (-8.2,1.7);
\node[red] at (-5.85,0.35) {$\frac{1}{\Delta_\omega}$};
\node at (-6.6,0.6) {$v$};
\node at (-8.6,1.9) {$w$};
\node at (-7.75,1.4) {$-$};
\draw  (-5.3,2.1) rectangle (-4.25,1.3);
\node at (-4.8,1.7) {$G$};
\node at (-6.8,1.7) {$C$};
\draw [-triangle 60](-3.45,1.7) -- (-2.65,1.7);
\node at (-3.15,2) {$y$};
\draw  (-3.65,1.7) ellipse (0.2 and 0.2);
\draw [-triangle 60](-6.3,1.7) -- (-5.3,1.7);
\draw [-triangle 60](-7.8,1.7) -- (-7.3,1.7);
\draw [-triangle 60](-4.25,1.7) -- (-3.85,1.7);
\draw [-triangle 60](-3.65,2.4) -- (-3.65,1.9);
\node at (-3.85,2.4) {$d$};
\end{tikzpicture}
	\caption{Feedback interconnection with feedback damper}
	\label{fig:schema_inv}
\end{figure}

The scheme of the second application, FMC, is in Fig.~\ref{fig:schema_inv}. The proposed concept adapts an {\em inverse shaper} application elaborated in \cite{Vyhlidal2016Feedback}. A typical application of this concept is position-control of a crane trolley ($G$: main body) with the aim to compensate the oscillatory modes of the suspended payload ($F$: flexible subsystem), i.e., the payload should not sway once the main body position $y$ reaches the set-point value $w$. The architecture in Fig.~\ref{fig:schema_inv} ensures the mode compensation also in the responses to the main-body disturbance $d$. For the crane application the mode of $F(s)$ to be compensated is assumed $\pm i\omega$, where $\omega=\sqrt{\frac{g}{L}}$, $L$ is the length of the payload and $g$ is gravitational acceleration. 

The adaptation of the concept is in substituting the inverse shaper by the transfer function $\frac{1}{\Delta_\omega(s)}$. As can be seen from the transfer functions
\begin{equation}
    T_{y_sw}=\frac{C(s)G(s)\Delta_\omega(s)}{\Delta_\omega(s)+C(s)G(s)}F(s),
\end{equation}
\begin{equation}
    T_{y_sd}=\frac{\Delta_\omega(s)}{\Delta_\omega(s)+C(s)G(s)}F(s),
\end{equation}
with $C(s)$ denoting the feedback controller, the double root at $\pm i\omega$ compensates the oscillatory pole of $F(s)$. Analogously to the previous application, the root multiplicity two enhances the robustness in mode compensation. Let us note that if the mode to be compensated is damped, i.e. given by $-\zeta\omega\pm i \omega\sqrt{1-\zeta^2}$, the parameters of $\Delta_\omega(s)$ can be adapted according to \eqref{CoeffsComplex}. 

For both above potential applications, only the concept was outlined with the simplest possible structure of $\Delta(s)$ obtained for a purposefully selected delay value $\tau$. A more detailed analysis is needed, mainly in studying the stability posture/margin of the overall systems with respect to the delay length. Possibly, selection of delay satisfying $0<\tau<\frac{\pi}{\omega}$ can be beneficial. Then the parameter determining rules \eqref{CoeffsComplex} are needed in their full complexity.

\subsection{Equations of higher order}
\label{SecHigherOrder}

One may consider, instead of \eqref{DiffEqn}, a $n$-th order equation with derivatives of order up to $n-1$ in the delays, and the corresponding quasipolynomial $\Delta$ of degree $2n$ made of a $n$-th degree polynomial and a polynomial of degree $n-1$ multiplied by $e^{-s\tau}$. The problem of assigning a real root of multiplicity $2n$ and proving its dominance has already been considered in \cite{MazantiMultiplicity, MazantiQualitative}. As for the assignment of complex conjugate roots of multiplicity $n$ each and proving their dominance, several arguments used in the present paper still hold with only minor modifications. For instance, the proof of Theorem~\ref{MainTheo}\ref{DominanceComplex} only requires continuity of the coefficients \eqref{CoeffsComplex} of the quasipolynomial with respect to $\theta_0$ as well as the MID property for the case of a real root of multiplicity $2n$. The main difficulty in generalizing the results of this paper to equations of higher order relies on providing suitable characterizations of the coefficients. Explicit characterizations such as \eqref{CoeffsComplex} seem intractable in the general case, but one may still rely on implicit characterizations, such as those in \cite[Lemma~1]{Boussaada2016Tracking}.

\section{Concluding remarks}

We have considered in this paper the multiplicity-induced-dominancy property for the linear time-invariant delay differential equation \eqref{DiffEqn} when placing a pair of complex conjugate roots of maximal multiplicity of its characteristic quasipolynomial \eqref{Delta}. Our main result, Theorem~\ref{MainTheo}, provides necessary and sufficient conditions for a pair of complex conjugate numbers $\sigma_0 \pm i \theta_0$ being roots of maximal multiplicity of \eqref{Delta} and shows that, under these conditions, these roots are necessarily strictly dominant.

We have also presented, in Section~\ref{SecRootsFuncTheta}, how other roots of \eqref{Delta} behave as the chosen roots $\sigma_0 \pm i \theta_0$ move away from the imaginary axis. Section~\ref{SecApplications} has illustrated the utility of our main result by presenting two engineering applications in which the proposed design is useful. Finally, we have discussed in Section~\ref{SecHigherOrder} how the results of this paper can be generalized to higher-order equations.

\begin{ack}
This work is supported by a public grant overseen by the French National Research Agency (ANR) as part of the ``Investissement d'Avenir'' program, through the iCODE project funded by the IDEX Paris-Saclay, ANR-11-IDEX-0003-02. The authors also acknowledge the support of Institut Polytechnique des Sciences Avanc\'ees (IPSA). The fourth author acknowledges support from the ESIF, EU Operational Programme Research, Development and Education, and from the Center of Advanced Aerospace Technology (CZ.02.1.01/0.0/0.0/16\_019/0000826), Faculty of Mechanical Engineering, Czech Technical University in Prague
\end{ack}

\bibliography{main}

\end{document}